\documentclass[twoside,leqno,twocolumn]{article} 

\usepackage{ltexpprt}
\usepackage{amsfonts,amssymb,latexsym}
\newcommand{\eqref}[1]{\hbox{$($\ref{#1}$)$}}
\newcommand{\operatorname}[1]{\hbox{#1}}
\newcommand{\binom}[2]{{#1\choose #2}}

\def\cal{\mathcal}
\def\R{\mathbb{R}}
\def\Z{\mathbb{Z}}
\def\Q{\mathbb{Q}}
\def\Pr{\mathbb{P}}
\def\Ex{\mathbb{E}}
\newcommand{\parag}[1]{\smallbreak{\bf\emph{#1}}}
\newcommand{\Img}[2]{\includegraphics[width=#1truecm]{#2}}

\def\Li{\operatorname{Li}}
\def\ds{\displaystyle}
\def\la{\langle}\def\ra{\rangle}

\def\iff{\hbox{\bf if}}
\def\thenn{\hbox{\bf then}}
\def\elsee{\hbox{\bf else}}
\def\doo{\hbox{\bf do}}

\def\return{\hbox{return}}
\def\brk{\hbox{\bf break}}

\def\opn{\operatorname}

\def\gber{\Gamma\opn{B}}
\def\ggeo{\Gamma\opn{G}}
\def\gpoi{\Gamma\opn{P}}
\def\glog{\Gamma\opn{L}}
\def\gvn{\Gamma\operatorname{VN}[\cal P](\lambda)}
\newcommand{\GVN}[1]{\Gamma\operatorname{VN}[\cal #1](\lambda)}

\def\geo{\operatorname{Geo}}
\def\poi{\operatorname{Poi}}
\def\trie{\operatorname{trie}}

\def\U{{\bf U}}
\def\wh{\widehat}

\def\egf{{\sc egf}}
\def\ogf{{\sc ogf}}

\def\up{\upsilon}

\newcounter{noteno}
\newenvironment{note}
        {\refstepcounter{noteno}%
        \setlength{\baselineskip}{0.9\baselineskip}%
        \smallbreak\noindent%
        \begin{small}%
        {\bf Note~\arabic{noteno}.}}%
        {\par\smallbreak\end{small}}

\def\type{\operatorname{type}}

\newtheorem{definition}{Definition}

\usepackage{graphicx}

\begin{document}

\begin{twocolumn}

\title{\Large\bf On Buffon Machines and Numbers}

\author{Philippe Flajolet\thanks{{\sc Algorithms}, 
	INRIA, F-78153 Le Chesnay, France}
\and Maryse Pelletier\thanks{LIP6, UPMC, 4, Place Jussieu, 75005 Paris, France}
 \and Mich\`ele Soria${}^{\dagger}$}

\date{\small {\tt philippe.flajolet@inria.fr}, {\tt maryse.pelletier@lip6.fr},
{\tt michele.soria@lip6.fr}}


\maketitle

\begin{abstract}
\begin{small}
The well-know   needle experiment  of Buffon 
can be regarded as an analog (i.e., continuous) device that 
stochastically ``computes'' the number
$2/\pi\doteq0.63661$, which is the experiment's probability of success.
Generalizing the  experiment and simplifying the computational framework, 
we consider probability  distributions,
which can be  produced \emph{perfectly},
from a  \emph{discrete source}  of un\-biased coin   flips.  We  describe  and
analyse  a few simple  \emph{Buffon   machines} that generate  geometric,
Poisson, and  logarithmic-series distributions.
We provide human-accessible Buffon  machines, which require a  dozen
coin flips  or   less,  on  average, and  produce    experiments whose
probabilities of success are expressible in terms of numbers  
such as~$\pi$, $\exp(-1)$, $\log2$,
$\sqrt{3}$, $\cos(\frac14)$, $\zeta(5)$. 
Generally, we develop a collection of constructions based on 
\emph{simple probabilistic mechanisms} that enable one to design
Buffon experiments involving 
compositions of exponentials and logarithms, polylogarithms, direct and inverse
trigonometric functions, algebraic and hypergeometric functions,
as well as functions defined by integrals, such as the Gaussian error function.
\par
\end{small}
\end{abstract}


%

%

\section*{Introduction}
Buffon's experiment (published in 1777)
is as follows~\cite{Badger94,Buffon77}.
\emph{Take a plane marked with parallel lines at unit distance from one another; 
throw a needle at random; finally, declare the experiment a success
if the needle intersects one of the lines.}
Basic calculus implies that the probability of success is $2/\pi$.
%

One can regard Buffon's experiment as a simple analog device that takes as input 
\emph{real} uniform $[0,1]$--random variables 
(giving the position of the centre and the angle of the needle)
and  outputs a discrete $\{0,1\}$--random variable,
with~$1$ for success and~$0$ for failure.
The process then involves \emph{exact} 
arithmetic operations over real numbers.
In the same vein,
the classical problem of \emph{simulating} random variables can be described
as the construction of analog devices (algorithms)
operating over the \emph{reals} and equipped
with a \emph{real} source of randomness,
which are both simple and computationally efficient.
The encyclopedic treatise  of Devroye~\cite{Devroye86b} 
provides many examples relative to 
the simulation of distributions, such as
Gaussian, exponential, Cauchy, stable, and so on.

\begin{figure}[t]\small
\begin{center}\renewcommand{\tabcolsep}{2truept}\renewcommand{\arraystretch}{1.5}
\begin{tabular}{ll|c|c|c}
\hline\hline
\multicolumn{2}{c}{\emph{Law}}
 & \emph{supp.} & \emph{distribution} 
& \emph{gen.} \\
\hline
 Bernoulli & Ber$(p)$  &$\{0,1\}$ & $\Pr(X=1)=p$ 
&$\gber(p)$ \\
 geometric & Geo$(\lambda)$ &$\Z_{\ge0}$ & $\Pr(X=r)=\lambda^r(1-\lambda)$ 
& $\ggeo(\lambda)$ \\
Poisson,& Poi$(\lambda)$ & $\Z_{\ge0}$  & $\ds \Pr(X=r)=e^{-\lambda}\frac{\lambda^r}{r!}$   
& $\gpoi(\lambda)$ \\[1.5mm]
logarithmic, & Log$(\lambda)$ & $\Z_{\ge0}$ & $\ds \Pr(X=r)=\frac{1}{L_{\vphantom{I_I}}}
\frac{\lambda^r}{r}$   
& $\glog(\lambda)$  \\
\hline\hline
\end{tabular}
\end{center}
\vspace*{-2.5truemm}
\caption{\label{tab-fig}\small
Main discrete probability laws: support, expression of the probabilities,
and naming convention for generators ($L:=\log(1-\lambda)^{-1}$).}
\end{figure}

\parag{Buffon machines.}
Our objective is   the \emph{perfect simulation of discrete random
variables} (i.e., variables supported by  $\Z$ or one of its  subsets;
see Fig.~\ref{tab-fig}). 
In this context, it is natural to start with a discrete source of randomness that
produces  \emph{uniform random bits} (rather than uniform $[0,1]$ real numbers),
since we are
interested  in   finite  computation  (rather than  infinite-precision
real numbers); cf Fig.~\ref{bmach-fig}.

\begin{definition}\label{bm-def}
A \emph{Buffon machine} is a 
deterministic device belonging to a 
computationally universal class (Turing machines, equivalently, register machines),
equipped with an external 
\emph{source of independent uniform random bits}
and input--output queues capable of storing integers
(usually, $\{0,1\}$-bits), 
which is assumed to halt\footnote{Machines that
\emph{always} halt can only produce Bernoulli distributions whose
parameter is a dyadic rational~$s/2^t$; see~\cite{KnYa76} and~\S\ref{frame-sec}.} \emph{with probability~$1$}.
\end{definition}
The fundamental question is then the following.
\emph{How does one generate, \emph{exactly} and \emph{efficiently}, 
discrete distributions 
using only  a discrete source of random  bits and finitary computations?}
A pioneering work in this direction is that of Knuth and Yao~\cite{KnYa76}
who discuss the power of various restricted devices (e.g., finite-state machines).
Knuth's treatment~\cite[\S3.4]{Knuth98}
and the articles~\cite{FlSa86,Yao85} provide additional results along these lines.

\begin{figure}\small
\begin{center}
\setlength{\unitlength}{1truecm}
\begin{picture}(6.5,5.5)
\put(0,0){\Img{6.5}{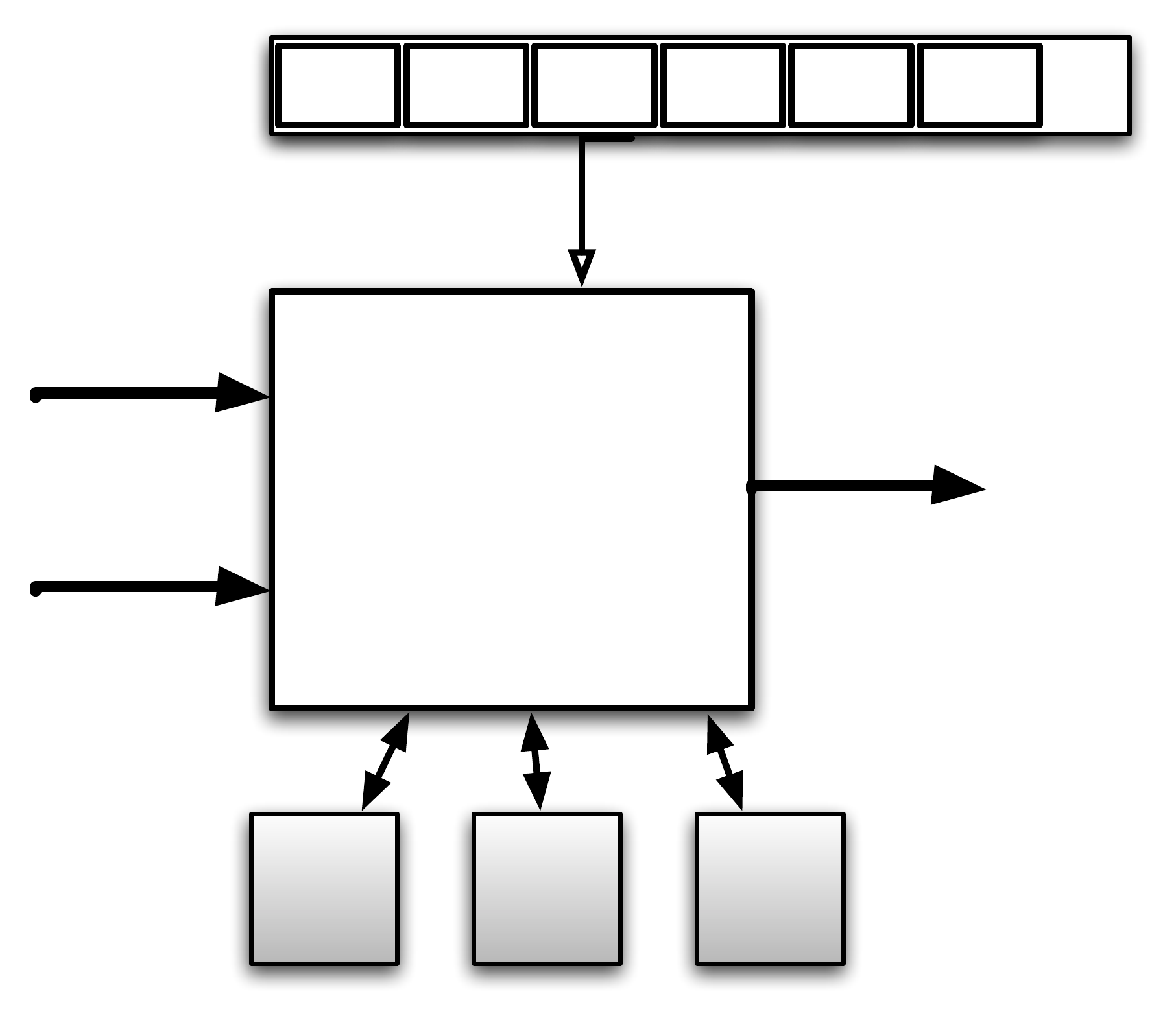}}
\put(-0.5,2.8){{\tt inputs}}
\put(5.75,2.8){{\tt output}}
\put(0.2,5.1){\tt flips}
\put(5.1,0.6){\tt Registers}
\put(1.75,5.05){\bf H}
\put(2.45,5.05){\bf T}
\put(3.15,5.05){\bf T}
\put(2.5,2.8){\begin{minipage}{1.5cm}{\em Finite\\control}\end{minipage}}
\end{picture}
\end{center}
\vspace*{-7truemm}
\caption{\small\label{bmach-fig} A Buffon machine with two inputs,
one output, and three registers.}
\end{figure}

Our original motivation for studying discrete Buffon machines came from
\emph{Boltzmann samplers} for combinatorial structures,
following the approach of Duchon, Flajolet, Louchard, and Schaeffer~\cite{DuFlLoSc04}
and Flajolet, Fusy, and Pivoteau~\cite{FlFuPi07}. The current implementations
relie on real number computations, and they require
generating distributions such as geometric, Poisson, or logarithmic,
with various ranges of parameters---since the objects ultimately produced 
are discrete,
it is natural  to try and produce them by purely discrete means.

\parag{Numbers.}
Here is an intriguing range of related issues.
Let~$\cal M$ be an input-free  Buffon machine that outputs a random variable~$X$,
whose value lies in~$\{0,1\}$. 
It can be seen from the definition that such a machine 
$\cal M$, when called repeatedly,  produces an independent
sequence of \emph{Bernoulli random variables}.
We say that $\cal M$ is a 
\emph {Buffon machine} or \emph {Buffon computer} for the number $p:=\Pr(X=1)$.
We seek
\emph{simple mechanisms---Buffon machines---that produce,
from uniform $\{0,1\}$-bits,  Bernoulli variables whose probabilities of
success are numbers such as}
\begin{equation}\label{numbers}
1/\sqrt{2}, ~~ e^{-1}, ~~ \log 2, ~~ \frac{1}{\pi},
~~ \pi-3, ~~ \frac{1}{e-1}\,.
\end{equation}
This problem can be seen as a vast generalization of Buffon's needle problem,
adapted to the discrete world.

\parag{Complexity and simplicity.} 
We will impose 
the loosely defined constraint that the Buffon machines 
we consider be \emph{short} and
conceptually \emph{simple}, to the extent of being easily implemented by a human.
Thus, emulating infinite-precision computation with multiprecision interval arithmetics
or appealing to functions of high complexity as primitives is disallowed\footnote{
	The informal requirement of \emph{``simplicity''}
	can be captured by the formal 
	notion of \emph{program size}. All the programs we develop necessitate at
	most a few dozen register-machine instructions, see~\S\ref{concl-sec}
	and the Appendix, 
	as opposed to programs based on arbitrary-precision arithmetics,
	which must be rather huge; cf~\cite{ScGrVe94}.
	If program size  is unbounded, the
	problem becomes trivial, since
	any Turing-computable number~$\alpha$ can be emulated  
	by a Buffon machine 
	with a number of coin flips that is~$O(1)$ on average, 
	e.g., by computing the sequence
	of digits of~$\alpha$ on demand; see~\cite{KnYa76} 
	and Eq.~\eqref{bern0} below.
}.
Our Buffon programs only make use of simple integer counters,  
string registers (\S\ref{vn-sec}) and stacks (\S\ref{alg-sec}),
as well as ``bags'' (\S\ref{int-sec}). 
The reader may try her hand at determining 
Buffon-ness in this sense of some of the numbers listed in~\eqref{numbers}.
We shall, for instance,   produce  a Buffon computer 
for the constant~$\Li_3(1/2)$ of Eq.~\eqref{li3}
(which involves $\log 2$, $\pi$, and~$\zeta(3)$),
one that is human-compatible, that consumes on average less than 6~coin flips,
and that requires at most 20~coin flips in  95\% of the simulations.
We shall also devise \emph{nine different ways of simulating 
Bernoulli distributions whose probability involves~$\pi$,
some requiring as little as five coin flips on average}.
Furthermore, the constructions of this paper
can  all be encapsulated into a universal interpreter, 
briefly discussed in Section~\ref{concl-sec}, which 
 has less than 60 lines of
{\sc Maple} code and produces all the constants of Eq.~\eqref{numbers},
as well as many more.

In this extended abstract, we
focus the discussion on  \emph{algorithmic  design}. Analytic estimates can   be
approached by means of (probability, counting) generating functions in
the style of methods of  analytic combinatorics~\cite{FlSe09}; see
 the typical discussion in Subsection~\ref{vn1-subsec}. The main results 
are Theorem~\ref{poilog-thm} (Poisson and logarithmic generators),
Theorem~\ref{explog-thm} (realizability of exps, logs, and trigs),
Theorem~\ref{int-thm} (general integrator), 
and Theorem~\ref{invtrig-thm} (inverse trig functions).

\section{\bf Framework and examples} \label{frame-sec}

Our approach consists in setting up a system based on the 
\emph{composition} of simple probabilistic
experiments, corresponding to simple computing devices. 
(For this reason, the Buffon machines of Definition~\ref{bm-def}
are allowed input/output registers.)
The \emph{unbiased 
random-bit generator}\footnote{%
	This convention entails no loss of generality. Indeed, as first observed by von Neumann,
	suppose we only have a biased coin, where $\Pr(1)=p$, $\Pr(0)=1-p$, with $p\in(0,1)$. Then, 
	one should toss the coin twice: if {\tt 01} is observed, then output~0; if
	{\tt 10} is observed, then output~1; otherwise, 
repeat with new coin tosses. See, e.g., \cite{KnYa76,Peres92} for more on this topic.
}, with which Buffon machines are equipped
will be named ``{\bf flip}''. 
The \emph{cost measure} of a computation (simulation)
is taken to be the number of  flips.
(For the restricted devices we shall consider, 
the overall simulation cost is essentially proportional to this measure.)

\begin{definition}\label{real-def} 
The function
$\lambda\mapsto\phi(\lambda)$,  defined for $\lambda\in(0,1)$ and with
values $\phi(\lambda)\in(0,1)$, is   said to be \emph{weakly realizable} 
if  there  is  a
machine~$\cal M$, which,  when provided on its input with 
a perfect generator
of      Bernoulli     variables   of   
(unknown) parameter
$\lambda$,    outputs, with probability~$1$, a   Bernoulli    random variable  of
parameter~$\phi(\lambda)$. 
The function~$\phi$
 is said to be \emph{realizable}, resp., \emph{strongly realizable},
if the (random) number~$C$ of coin flips has finite expectation, 
resp., exponential tails\footnote{
$C\equiv C(\lambda)$ has exponential tails if
 there are  constants $K$ and $\rho<1$ such that $\Pr(C>m)\le K\rho^{m}$.}.
\end{definition}
\noindent 
We shall  also say  that  $\phi$  has a  [weak,
strong] \emph{simulation} if  it is realizable by  a  machine [in  the weak,
strong sense].
Schematically:
\[
\hbox{\setlength{\unitlength}{1truecm}
\begin{picture}(6,2.5)
\put(0.5,0){\Img{5.6}{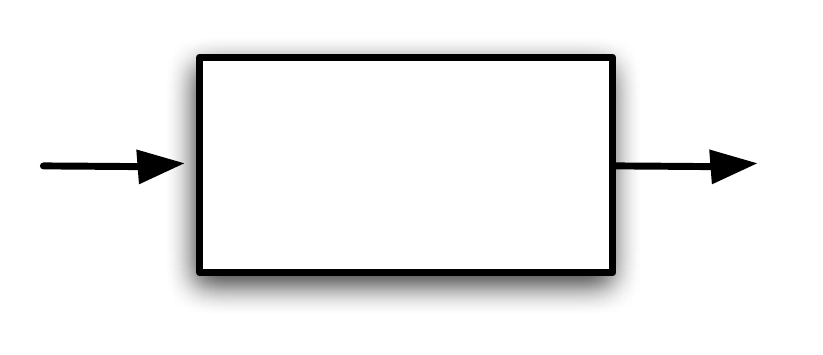}}
\put(-1.4,1.6){$X_1,X_2,\ldots\in\operatorname{Ber}(\lambda)$}
\put(5.0,1.6){$Y\in\operatorname{Ber}(\phi(\lambda))$}
\put(2.9,1.25){\large$\cal M$}
\end{picture}}
\]

\vspace*{-3truemm}

\noindent
The devices are normally implemented as \emph{programs}. 
Using $\gber(\lambda)$ as a generic notation for a Bernoulli generator of
parameter~$\lambda$, a Buffon machine that realizes the function~$\phi$
is then equivalent to a program
that can call (possibly several times) a $\gber(\lambda)$, as an external routine,
and then outputs a random Bernoulli variable of parameters~$\phi(\lambda)$.
It corresponds to  a construction of type
$\gber(\lambda)\longrightarrow \gber(\phi(\lambda))$.

The definition is extended to machines with $m$ inputs,
in which case a function $\phi(\lambda_1,\ldots,\lambda_m)$ 
of~$m$ arguments is \emph{realized}; see below for several examples.
A machine with no input register 
then computes  a  function  $\phi()$  of  no   argument,  that is,   a
constant~$p$, and it does so
based solely on its  source of unbiased coin flips: this
is  what  we called in  the  previous section  a  Buffon  machine  (or
computer) for~$p$. 

The fact that Buffon machines are allowed 
input registers makes it possible to \emph{compose}
them. For instance, if $\cal M$ and $\cal N$ realize the
unary functions $\phi$ and $\psi$, connecting the output of
$\cal M$ to the input of $\cal N$ realizes the 
composition $\psi\circ\phi$. It is one of our goals to
devise Buffon computers for special values of the success probability~$p$
by composition of simpler functions, eventually only applied 
to the basic $\operatorname{flip}$ primitive.

\smallskip

Note    that there are   obstacles to   what can   be done:  Keane and
O'Brien~\cite{KeOB94},   for  instance,   showed   that  the ``success
doubling''  function   $\min(2\lambda,1)$   cannot  be   realized  and
discontinuous  functions are    impossible   to simulate---this,   for
essentially measure-theoretic reasons. On the  positive side, Nacu and
Peres~\cite{NaPe05} show that  every (computable) Lipschitz
function is  realizable and  every (computable) real-analytic function
has   a  strong     simulation,  but  their      constructions require
\emph{computationally    unrestricted}   devices; namely,
sequences  of approximation 
functions   of  increasing   complexity.
Instead, our  purpose here is  to show how sophisticated \emph{perfect}
simulation   algorithms   can    be  systematically   and  efficiently
synthetized by
\emph{composition of \emph{simple} probabilistic processes}, this without 
the need of 
hard-to-compute sequences of approximations. 

To set the stage of the present study,
we shall briefly review in sequence:
$(i)$~decision trees and polynomial functions;
$(ii)$~Markov chains (finite graphs) and rational functions.

\begin{figure*}
\begin{equation}\label{polys}
\hbox{\small\begin{tabular}{lll}
\hline\hline
\emph{Name}  & \emph{realization} & \emph{function}\\
\hline
Conjunction ($P\wedge Q$) & \iff~$P()=1$ \thenn~\return($Q()$) \elsee~return(0) & $p\wedge q=p\cdot q$ \\
Disjunction ($P\vee Q$)  & \iff~$P()=0$ \thenn~\return($Q()$) \elsee~return(1) & $p \vee q = p+q-pq$ \\
Complementation ($\neg P$) & \iff~$P()=0$ \thenn~return(1) \elsee~return(0) & $1-p$\\
Squaring & ($P \wedge P$) & $p^2$\\
Conditional ($R\rightarrow P\,|\,Q$) & \iff~$ R()=1$ \thenn~return($P()$) \elsee~return($Q()$) & $rp+(1-r)q$.\\
Mean & \iff~flip() \thenn~return($P()$) \elsee~return($Q()$) & $\frac12(p+q)$. \\
\hline\hline
\end{tabular}}
\end{equation}
\hrule
\end{figure*}

\parag{Decision trees and polynomials.} 
Given three machines $P,Q,R$ with outputs in~$\{0,1\}$ and with $\Pr_P(1)=p$, $\Pr_Q(1)=q$, $\Pr_R(1)=r$,
we can easily build machines, corresponding to loopless programs,
whose probability of success is $p\cdot q$, $1-p$,
or any composition thereof: see Eq.~\eqref{polys} [top of next page]
 for the most important boolean primitives.
We can then simulate a Bernoulli of parameter any dyadic rational $s/2^t$, starting from the unbiased flip procedure,
performing~$t$ draws, and declaring a success in~$s$ designated cases.
Also, by calling a Bernoulli variate of (unknown) parameter~$\lambda$  a fixed number~$m$ of times, then observing
the sequence of outcomes and its number~$k$ of~$1$s,  we can realize any polynomial function 
$a_{m,k}\lambda^k(1-\lambda)^{m-k}$,
with $a_{m,k}$ any integer satisfying $0\le a_{m,k}\le \binom{m}{k}$.

\parag{Finite graphs (Markov chains) and rational functions.}
We now consider programs that allow iteration and correspond to a finite control graph---these are 
equivalent to Markov chains, possibly parameterized 
by ``boxes'' representing calls to
external Bernoulli generators. In this way, we can produce a Bernoulli  generator
of any rational parameter $\lambda\in\Q$, by combining a simpler dyadic generator with  iteration.
(For instance, to get a 
Bernoulli generator of parameter~$\frac13$, denoted by
$\gber(\frac13)$, flip an unbiased coin twice;
return success if {\tt 11} is observed, failure in case of {\tt 01} or {\tt 10}, and repeat the experiment in
case of {\tt 00}.) 
Clearly, only rational 
numbers and functions can be realized by finite graphs.

A highly useful construction in this range of methods is the \emph{even-parity}
construction:
%
\begin{equation}\label{evenp}
\hbox{\qquad}\hbox{\small\begin{minipage}{6.6truecm}
[even-parity]\quad \doo~$\{$\, \iff~$P()=0$ \thenn~return(1); 
\\ \hbox{\qquad\qquad\qquad\quad} \iff~$P()=0$ \thenn~return(0) $\}$.\end{minipage}}
\end{equation}

\noindent
This realizes the function $p\mapsto1/(1+p)$. Indeed,
the probability of~$k+1$ calls to~$P()$ is $(1-p)p^{k}$, which, when summed over $k=0,2,4,\ldots$,
gives the probability of success as $1/(1+p)$; thus, this function is realizable. Combining this function with complementation ($x\mapsto 1-x$)
leads, for instance, to a way of compiling any rational probability ${a}/{b}$
into a generator $\gber\left(\frac a b\right)$ whose size is
proportional to the sum of digits of the continued fraction representation
of $a/ b$. (See also Eq.~\eqref{bern0} below for an alternative based on 
binary representations and~\cite{KnYa76} for more on this topic.)

Here are furthermore two important characterizations based on works of Nacu--Peres~\cite{NaPe05} 
(see also W\"astlund~\cite{Wastlund99}) and Mossel--Peres~\cite{MoPe05}:

\begin{theorem}[\cite{MoPe05,NaPe05,Wastlund99}]
$(i)$~Any polynomial~$f(x)$ with rational coefficients that maps $(0,1)$ into~$(0,1)$ 
is strongly realizable by a finite graph.
$(ii)$~Any rational function~$f(x)$ with rational coefficients that maps $(0,1)$ into~$(0,1)$ 
is strongly realizable by a finite graph.
\end{theorem}
\noindent
(Part~$(i)$ is based on a theorem of P\'olya, relative to nonnegative polynomials;
Part~$(ii)$ depends on an ingenious ``block simulation'' principle.
The corresponding constructions however require unbounded precision arithmetics.)

We remark  that we can also produce a \emph{geometric variate} from 
a Bernoulli of the same parameter: just repeat till failure. This gives rise to the program
%
\begin{equation}\label{geom}
\hbox{\quad\quad}\hbox{\small
\begin{minipage}{7.3truecm}
 [Geometric]\quad $\ggeo(\lambda)$ := $\{$\,$K:=0$; \doo~$\{$\,\iff~$\gber(\lambda)=0$ 
\\ \hbox{\qquad\qquad\qquad\qquad}\thenn~return$(K)$; $K:=K+1$;\,$\}$\,$\}$.\end{minipage}
}
\end{equation}

\noindent
(The even-parity construction implicitly makes use of this.) 
The special  $\ggeo(\frac12)$ then simply iterates on the basis of a flip.

In case we have access to the complete binary expansion of~$p$
(note that
this is \emph{not} generally 
permitted in our framework), a Bernoulli generator is
obtained by
\begin{equation}\label{bern0}
\hbox{\qquad}\hbox{\small\begin{minipage}{7.25truecm}
[Bernoulli/binary]\quad 
$\{$\,{\bf let} $Z:=1+\ggeo(\frac12)$;
\\ \hbox{\qquad\qquad\qquad\qquad\qquad\quad
 return(bit$(Z,p)$)\,$\}$.}
\end{minipage}}
\end{equation}

\noindent
In words: \emph{in order to draw a Bernoulli variable of parameter~$p$
whose binary representation is available, return the 
bit of~$p$ whose random index is given by  a shifted geometric variable of parameter~$1/2$.}
(Proof: emulate a comparison 
between a uniformly random $V\in[0,1]$ with~$p\in[0,1]$; see~\cite[p.~365]{KnYa76}
for this trick.)
The cost is by design  a geometric of rate $1/2$.
In particular, Eq.~\eqref{bern0}
automatically gives us  a $\gber(p)$, for any rational~$p\in\Q$, 
by means of a simple Markov chain,
based on the eventual periodicity of the representation of~$p$. 
(The construction will prove useful when we discuss ``bags'' in \S\ref{int-sec}.)

%
%
%

\section{\bf The von Neumann schema}\label{vn-sec}

The idea of generating  certain probability distributions by  way of
their Taylor expansion  seems to  go back  to von  Neumann.  An  early
application discussed  by   Knuth and  Yao~\cite{KnYa76},  followed by
Flajolet and Saheb~\cite{FlSa86},    is  the \emph{exact}  simulation  of  an
exponential variate by means of random $[0,1]$--uniform  variates, this by a ``continuation'' process
that altogether \emph{avoids multiprecision operations}. 
We build upon von Neumann's idea and introduce 
in Subsection~\ref{vn1-subsec}
a general schema for random generation---the \emph{von Neumann schema}.
We then explain in Subsection~\ref{vn2-subsec}
how this schema may be adapted to 
realize classical transcendental functions, such as $e^{-\lambda}$, $\cos(\lambda)$,
only knowing a generator~$\gber(\lambda)$.

\subsection{Von Neumann generators of discrete distributions.} \label{vn1-subsec}
First a few notations from~\cite{FlSe09}. 
Start from a class~$\cal P$ of permutations, with~$\cal P_n$
the subset of permutations of size~$n$ and $P_n$
the cardinality of $\cal P_n$. The 
(counting) \emph{exponential generating function}, or {\egf}, is
\[
P(z):=\sum_{n\ge 0} P_n \frac{z^n}{n!}.
\]
For instance, the classes $\cal Q, \cal R,\cal S$ of, respectively, 
\emph{all} permutations, 
\emph{sorted} permutations, and \emph{cyclic} permutations have {\egf}s given by
$
Q(z)=(1-z)^{-1},~
R(z)=e^z, ~
S(z)=\log(1-z)^{-1}$,
since $Q_n=n!$, $R_n=1$, $S_n=(n-1)!$, for $n\ge1$.
We observe that the class~$\cal S$ of cyclic permutations is isomorphic to
the class of permutations such that the maximum occurs in the first position:
$U_1>U_2,\ldots,U_N$. (It suffices to ``break'' a cycle at its maximum element.)
We shall also denote by~$\cal S$ the latter class, which is easier 
to handle algorithmically.

Let $\U=(U_1,\ldots,U_n)$ be a vector of real numbers. 
By replacing each $U_j$ by its rank in~$\U$, we obtain a permutation
$\sigma=(\sigma_1,\ldots,\sigma_n)$, which is called the (order) \emph{type} of 
$\U$ and is written $\type(\U)$.
For instance: $\type(1.41,0.57,3.14,2.71)=(2,1,4,3)$. 
The von Neumann schema, relative to a class $\cal P$ of permutations is
described in Fig.~\ref{VN-fig} and denoted by~$\gvn$. Observe that it only needs 
a geometric generator $\ggeo(\lambda)$, hence, eventually,
 only a Bernoulli generator $\gber(\lambda)$, whose parameter~$\lambda$ is 
\emph{not} assumed to be known.

\begin{figure}\small

\begin{tabbing}
\=X \= X \= XXX\= XXX \= XXX\=\kill
\>\>$\gvn$ := $\{$\, \doo~$\{$ 
\\
\>\>\>\> $N:=\ggeo(\lambda)$; 
\\
\>\>\>\> {\bf let} $\U:=(U_1,\ldots,U_N)$, vector of $[0,1]$--uniform 
\\
\>\>\>\> \emph{$\{$\,bits of the $U_j$ are produced on  call-by-need basis$\}$}
\\
\>\>\>\> {\bf let} $\tau:=\trie(\U)$; 
{\bf let} $\sigma:=\type(\U)$; \\
\>\>\>\> \iff~$\sigma\in\cal P_N$ \thenn~return($N$)\,$\}$\,$\}$.
\end{tabbing}

\vspace*{-4.5truemm}
\caption{\label{VN-fig}\small
The von Neumann schema $\gvn$, in its basic version,
relative to a class of permutations~$\cal P$ and a 
parameter~$\lambda$ (equivalently given by its Bernoulli generator~$\gber(\lambda)$).}
\end{figure}

First,  by construction, a value~$N$ is,
at each stage of the iteration, chosen with probability~$(1-\lambda)\lambda^N$. 
The procedure consists in a sequence of failed trials (when
$\type({\bf U})$ is not in~$\cal P$), followed by eventual success.
An iteration (trial) 
that succeeds then returns the value $N=n$ with probability
\begin{equation}\label{gf0}
\frac{(1-\lambda)P_n\lambda^n/n!}{(1-\lambda)\sum_n P_n\lambda^n/n!}
=\frac{1}{P(\lambda)}\frac{P_n\lambda^n}{n!}.
\end{equation}
For $\cal P$ one of the three classes $\cal Q,\cal R,\cal S$ described 
above this gives us three interesting distributions:
\begin{equation}\label{gf1}
\hbox{\small\begin{tabular}{c|c|c}
\hline\hline
  all $(\cal Q)$ & sorted $(\cal R)$& cyclic $(\cal S)$\\
\hline
 $(1-\lambda)\lambda^n$ &
$\ds e^{-\lambda}\frac{\lambda^{n^{\vphantom{X}}}}{n!}$ & $\ds \frac{1}{L}\frac{\lambda^n}{n}$
\\
geometric & Poisson & logarithmic. \\
\hline\hline
\end{tabular}}
\end{equation}
%
%
%
\noindent
The case of all permutations ($\cal Q$) is trivial, since  
no order-type  restriction is involved, so that the initial
value of $N\in\geo(\lambda)$ is returned. It is  notable 
  that, in the other two cases $(\cal R,\cal S)$, \emph{one produces
the Poisson and logarithmic distributions,
by means of permutations obeying simple restrictions}.

Next, implementation details should be discussed. 
Once~$N$ has been drawn,
we can imagine producing the~$U_j$ in sequence, by generating at each stage only the \emph{minimal}
number of bits needed to distinguish any $U_j$ from the other ones.
This corresponds to the construction of a \emph{digital tree},
also known as \emph{``trie''}~\cite{Knuth98a,Mahmoud92,Szpankowski01}
and is summarized by the command ``{\bf let} $\tau:=\trie(U)$'' in
the schema of Fig.~\ref{VN-fig}. As the digital tree~$\tau$ is constructed,
the~$U_j$ are gradually sorted, so that the order type~$\sigma$ can be 
determined---this involves no additional flips, just bookkeeping.
The test $\sigma\in\cal P_N$ is of this nature and it requires no flip at all.

The general properties of the von Neumann schema are summarized
as follows.

\begin{theorem} \label{vn-thm}
$(i)$~Given an arbitrary
 class~$\cal P$ of permutations and a parameter~$\lambda\in(0,1)$,
the von Neumann schema $\gvn$ produces \emph{exactly} a 
discrete random variable with probability distribution
\[
\Pr(N=n)=\frac{1}{P(\lambda)}\frac{P_n\lambda^n}{n!}.
\]

$(ii)$~The number~$K$ of iterations has expectation $1/s$, where
$ 
s={(1-\lambda)P(\lambda)},
$ 
and its distribution is $1+\geo(1-s)$. 

$(iii)$~The number~$C$ of flips consumed by the algorithm
(not counting
the ones in $\ggeo(\lambda)$)
is a random variable with probability generating function
\begin{equation}\label{didi}
\Ex(q^C)=\frac{H^+(\lambda,q)}{1- H^-(\lambda,q)},
\end{equation}
where $H^+,H^-$ are computable from the family
of polynomials in~\eqref{recu} below by
\[\begin{array}{lll}
H^+(z,q)&=&\ds (1-z)\sum_{n=0}^\infty \frac{P_n}{n!} h_n(q)z^n\\
H^-(z,q)&=&\ds (1-z)\sum_{n=0}^\infty  \left(1-\frac{P_n}{n!}\right) h_n(q)z^n.
\end{array}\]
The distribution  of~$C$ has exponential tails.
\end{theorem}
\begin{small}\begin{proof} [Sketch]
Items~$(i)$ and~$(ii)$ result from the discussion above.
Regarding Item~$(iii)$, a crucial observation is that the digital tree 
created at each step of the iteration is only a function of the
underlying \emph{set} of values. But there is complete independence between
this \emph{set} of values and their \emph{order type}. This 
justifies~\eqref{didi}, where $H^+,H^-$ are the probability generating functions associated with success and failure of one trial, respectively.

We next  need to discuss the fine structure of costs.
The cost of each iteration,
as measured by the
number of coin flips, is exactly that of generating the tree~$\tau$
of random size~$N$. The number of coin flips to build~$\tau$ 
coincides with the \emph{path length} of~$\tau$, written~$\omega(\tau)$,
which satisfies the
inductive definition
\begin{equation}\label{plindu}
\omega(\tau)=|\tau|+\omega(\tau_0)+\omega(\tau_1),
\qquad |\tau|\ge 2,
\end{equation}
where $\tau=\la \tau_0,\tau_1\ra$ and $|\tau|$ is the size of~$\tau$, that is,
the number of elements contained in~$\tau$.

Path length  is a much studied parameter, starting with
the work of Knuth in the mid 1960s relative to the analysis of 
radix-exchange sort~\cite[pp.~128--134]{Knuth98a}; see also
the books of Mahmoud~\cite{Mahmoud92} and Szpankowski~\cite{Szpankowski01}
as well as
 Vall\'ee \emph{et al.}'s analyses 
under dynamical source models~\cite{ClFlVa01,Vallee01}.
It is known
from these works that the expectation of path length, for~$n$ uniform  binary sequences,
is \emph{finite}, with exact value
\[
\Ex_{n}[\omega]=n\sum_{k=0}^\infty \left[1-\left(1-\frac{1}{2^k}\right)^{n-1}
\right],
\]
and asymptotic form (given by a Mellin transform analysis~\cite{FlGoDu95,Mahmoud92,Szpankowski01}):
$
\Ex_{n}[\omega]=n\log_2 n + O(n).
$
A consequence of this last fact is that $\Ex(C)$ is finite,
i.e., the generator~$\gvn$ has the basic
\emph{simulation (realizability) property}.

The distribution of path length is known to be asymptotically Gaussian, as~$n\to\infty$,
after Jacquet and R\'egnier~\cite{JaRe86}
and Mahmoud \emph{et al.}~\cite{MaFlJaRe00};
see also~\cite[\S11.2]{Mahmoud00}. 
For our purposes, it suffices to note that the bivariate {\egf} 
\[
H(z,q):=\sum_{n=0}^\infty \Ex_{n}[q^\omega] \frac{z^n}{n!}
\]
satisfies the nonlinear functional equation
$ 
H(z,q)=H\left(\frac{zq}{2},q\right)^2+z(1-q)$, with $ H(0,q)=1$.
This equation fully determines~$H$, since it is equivalent to a
recurrence on coefficients,
$h_n(q):=n![z^n]H(z,q)$,
for $n\ge2$:
\begin{equation}\label{recu}
\hbox{\quad}h_n(q)=\frac{1}{1-q^n 2^{1-n}}\sum_{k=1}^{n-1} 
\frac{1}{2^n}\binom{n}{k} h_k(q)h_{n-k}(q).
\end{equation}
The computability of~$H^+,H^-$ then results.
In addition, large deviation estimates  
can be deduced from~\eqref{recu}, which serve to establish exponential
tails for~$C$, thereby ensuring a \emph{strong simulation property}
in the sense of Definition~\ref{real-def}.
\end{proof}\par
\end{small}

A notable consequence of Theorem~\ref{vn-thm}  is the 
possibility of generating a Poisson or logarithmic variate 
by a simple device: as we saw
in the discussion preceding the statement of the Theorem,
only \emph{one branch} of the trie needs to be maintained,
in the case of the classes~$\cal R$ and~$\cal S$ of~\eqref{gf1}.

\begin{theorem}\label{poilog-thm} 
The Poisson and logarithmic distributions of parameter~$\lambda\in(0,1)$
have a strong simulation by a  Buffon  machine,
$\GVN{R}$ and~$\GVN{S}$, respectively, 
which only uses  a
single string register.
\end{theorem}

Since the sum of two Poisson variates is Poisson (with rate the sum of the rates),
the strong simulation result extends to  any $X\in\poi(\lambda)$, for any $\lambda\in\R_{\ge0}$.
\emph{This answers an open question of Knuth and Yao} in~\cite[p.~426]{KnYa76}.
We may also stress here that the distributions of costs are easily computable:
with the symbolic manipulation system {\sc Maple},
the cost of generating a Poisson(1/2) variate
is found to have probability generating function (Item~$(iii)$ of Theorem~\ref{vn-thm})
\[
{\frac {3}{4}}+{\frac {7}{128}}{q}^{2}+{\frac {119}{4096}}{q}^{4}+{
\frac {19}{1024}}{q}^{5}+{\frac {2023}{131072}}{q}^{6}+
\frac{179}{16384}q^7+
\cdots\,.
\]
Interestingly enough, the analysis of the logarithmic generators involves ideas similar to those
relative to a classical leader election protocol~\cite{FiMaSz96,Prodinger93}.

\subsection{Buffon computers: logarithms, exponentials, and trig functions.} \label{vn2-subsec}
We can also take any of the previous constructions and specialize 
it by declaring a success whenever a special value $N=a$ is returned,
for some predetermined $a$ (usually $a=0,1$),
declaring a failure, otherwise. For instance, the Poisson generator with~$a=0$
gives us in this way a Bernoulli generator with parameter
$
\lambda'=\exp(-\lambda)
$.
Since the von Neumann machine only requires a Bernoulli generator $\gber(\lambda)$,
we thus have a purely discrete construction 
$ 
\gber(\lambda) \longrightarrow \gber\left(e^{-\lambda}\right).
$ 
Similarly, the logarithmic generator 
restricted to~$a=1$ provides a construction
$ 
\gber(\lambda) \longrightarrow \gber\left(\frac{\lambda}{\log(1-\lambda)^{-1}}\right).
$ 
Naturally, these constructions can be enriched by the basic ones of Section~\ref{frame-sec},
in particular, complementation.

Another possibility is to make use of the number $K$ of iterations,
which is a shifted geometric of rate $s=(1-\lambda)P(\lambda)$;
see Theorem~\ref{vn-thm}, Item~$(ii)$.
If we declare a success when~$K=b$, for some predetermined~$b$, we then obtain 
yet another brand of generators. The Poisson generator 
used in this way with~$b=1$ gives us
$ 
\gber(\lambda)~ \longrightarrow~
\gber\left((1-\lambda)e^{\lambda}\right),~
\gber\left(\lambda e^{1-\lambda}\right),
$ 
where   the   latter involves   an additional   complementation.    

%

Trigonometric functions can also be brought into the game.
A sequence $\U=(U_1,\ldots,U_n)$ is said to be \emph{alternating}
if $U_1<U_2>U_3<U_4>\cdots$. 
It is well known that the {\egf}s of the classes $\cal A^+$ of even-sized and 
$\cal A^-$ of odd-sized permutations 
are respectively
$
A^+(z)=\sec(z)={1}/{\cos(z)}$,
and $
A^-(z)=\tan(z)={\sin(z)}/{\cos(z)}$.
(This result, due to D\'esir\'e Andr\'e around 1880, is in most books
on combinatorial enumeration, e.g., ~\cite{FlSe09,GoJa83,Stanley99}.)
Note that the property of being alternating can 
once more be tested sequentially: only 
a partial expansion of the current value of~$U_j$ needs to be stored at any given instant.
By making use of the properties $\cal A^+, \cal A^-$, with, respectively $N=0,1$, we 
then obtain new trigonometric constructions. 

In summary:

\begin{theorem}\label{explog-thm}
 The following functions admit a strong simulation:
\[
\begin{array}{l}
e^{-x},~~ e^{x-1},~~(1-x)e^x,~~xe^{1-x}, \\
\ds \frac{x}{\log(1-x)^{-1}},~~ \frac{1-x}{\log(1/x)},~~ (1-x)\log\frac{1}{1-x},~~
x\log\frac1x, \\[3.5truemm]
\ds {\cos(x)}, ~~\frac{1-x}{\cos(x)},~~\frac{x}{\tan(x)},~~ (1-x)\tan(x).
\end{array}
\]
\end{theorem}

%
%
%
%
%
%


\subsection{\bf Polylogarithmic constants.}
The probability that a vector~$\U$ is such that $U_1>U_2,\ldots,U_n$ (the first element
is largest)
equals $1/n$, a property that underlies the logarithmic series generator.
By sequentially drawing~$r$ several such vectors and requiring success on \emph{all}~$r$ trials,
we deduce constructions for families involving the polylogarithmic functions,
$
\Li_r(z):=\sum_{n\ge1}^\infty {z^n}/{n^r}$,
 with $ r\in \Z_{\ge1}$.
Of course, $\Li_1(1/2)=\log 2$.
The few special evaluations known for polylogarithmic values
 (see the books by Berndt on Ramanujan~\cite[Ch.~9]{Berndt85}
and by Lewin~~\cite{Lewin81,Lewin91})
include $\Li_2(1/2)=\frac{\pi^2}{12}-\frac{1}{2}\log^2 2$ and
\begin{equation}\label{li3}
\hbox{\hspace*{0.7truecm}\small
$\ds\Li_3(\frac{1}{2})=\frac{\log^3 2}{6}-\frac{\pi^2\log 2}{12}+\frac{7\zeta(3)}{8},
~ \zeta(s):=\sum_{n\ge1}\frac{1}{n^s}.$}
\end{equation}
By rational convex combinations, we obtain Buffon computers for
$
\frac{\pi^2}{24}$ and $ \frac{7}{32}\zeta(3).
$
Similarly, the celebrated BBP (Bailey--Borwein--Plouffe) formulae~\cite{BaBoPl97}
can be implemented as Buffon machines.

\section{\bf Square roots, algebraic, and hypergeometric functions} \label{alg-sec}
We now examine a new brand of generators
based on properties of \emph{ballot sequences},
which open the way to new constructions,
including an important square-root mechanism.
The probability that, in~$2n$  tosses of a fair coin, there are as many heads as tails is
$\varpi_n = \frac{1}{2^{2n}}\binom {2n}{n}$.
The property is easily testable with a single integer counter~$R$ subject
only to the basic operation $R:= R\pm 1$ and to the basic test $\ds R\mathop{=}^{?} 0$.
From this, one can build a square-root computer and, by repeating the test,
certain hypergeometric constants can be obtained.

%

\subsection{Square-roots}
Let~$N$ be  a random variable with distribution $\geo(\lambda)$. Assume we flip
$2N$~coins and return a success, if the score of heads and tails is balanced.
The probability of success is
\[
S(\lambda):=\sum_{n=0}^\infty (1-\lambda)\lambda^n \varpi_n =\sqrt{1-\lambda}.
\]
(The final simplification is due to the binomial expansion of $(1-x)^{-1/2}$.)
This simple property gives rise to the \emph{square-root construction}
due to W\"astlund~\cite{Wastlund99}
and  Mossel--Peres~\cite{MoPe05}:

\begin{small}\begin{tabbing}
X \= XXXX \= XXX \= XXX\= XXX \= XXX\=\kill
$\gber\left(\sqrt{1-\lambda}\right)$ := 
 \>\>\> $\{$\,{\bf let} $N:=\ggeo(\lambda)$; 
\\
\> {\bf draw} $X_1,\ldots, X_{2N}$ {with} $\Pr(X_j=+1)=\Pr(X_j=-1)=\frac12$;\\
\> {\bf set} $\Delta:=\sum_{j=0}^{2N} X_j$; 
\\ 
\> \iff~$\Delta=0$ \thenn~return(1) \elsee~return(0)\,$\}$.
\end{tabbing}
\end{small}

\noindent
The mean number of coin flips used is then simply obtained by differentiation
of generating functions.

\begin{theorem}[\cite{MoPe05,Wastlund99}]\label{sqrt-thm}
The square-root construction 
yields a Bernoulli generator of parameter~$\sqrt{1-\lambda}$,
given a $\gber(\lambda)$. The mean number of coin flips required, not counting the ones
involved in the calls to~$\gber(\lambda)$, is $
\frac{2\lambda}{1-\lambda}$. The function $\sqrt{1-\lambda}$ is strongly realizable.
\end{theorem}
\noindent
By complementation of the original Bernoulli generator, we also have a construction
$
\gber(\lambda) ~\longrightarrow~ \gber(1-\lambda)~
\longrightarrow \gber\left(\sqrt{\lambda}\right),
$
albeit one that is irregular at~$0$.

\begin{note} \emph{Computability with a pushdown automaton.} \label{pda-note}
It can be seen that the number~$N$ in the square-root generator never needs to be
stored explicitly: an equivalent form is
\begin{tabbing}
\= XXXX \= XXX \= XX\= XXX \= XXX\=\kill
\>$\gber\left(\sqrt{1-\lambda}\right)$ :=  \>\>\> $\{$\,\doo~$\{$\, $\Delta:=0$;\\
\>\> \iff~$\gber(\lambda)=0$ \thenn~\brk; \\
\>\> \iff~flip=1 \thenn~$\Delta:=\Delta+1$ \elsee~$\Delta:=\Delta-1$; 
\\
\>\> 
\iff~flip=1 \thenn~$\Delta:=\Delta+1$ \elsee~$\Delta:=\Delta-1$\,$\}$\\
\>\> \iff~$\Delta=0$ \thenn~return(1) \elsee~return(0)\,$\}$.
\end{tabbing}
In this way, only a stack of unary symbols needs to be maintained: the stack keeps track of the 
absolute value~$|\Delta|$ stored in unary, the finite control can keep track of the sign
of~$\Delta$. We thus have \emph{realizability of the square-root construction
by means of a pushdown (stack) automaton}.

This suggests a number of variants of the square-root construction,
which are also computable by a pushdown automaton. For instance,
assume that,
upon the condition ``flip=1'', one does $\Delta:=\Delta+2$ (and still does $\Delta:=\Delta-1$ otherwise).
The sequences of ~$H$ (heads) and $T$ (tails) that lead to
eventual success (i.e., the value~1 is returned) 
correspond to lattice paths that are bridges with vertical steps in $\{+2,-1\}$;
see~\cite[\S VII.8.1]{FlSe09}. The corresponding counting generating function
is then
$ 
S(z)=\sum_{n\ge0} \binom{3n}{n} z^{3n},
$ 
and the probability of success is
$
(1-\lambda)S\left(\frac{\lambda}{2}\right)
$.
As it is well known (via Lagrange inversion), the function~$S(z)$ is a variant of 
an algebraic function of degree~$3$; namely,
$ 
S(z)=(1-3zY(z)^2)^{-1}$, where $Y(z)=z+zY(z)^3,
$ 
and $Y(z)=z+z^4+3z^7+\cdots$ is a generating function of ternary trees. One can synthetize in
the same way the family of algebraic functions
\[
S(z)\equiv S^{[t]}(z)=\sum_{n\ge0} \binom{tn}{n} z^{tn},
\]
by updating ~$\Delta$ with $\Delta:=\Delta+(t-1)$. 
\end{note} 

As a consequence of Theorems~\ref{explog-thm} and~\ref{sqrt-thm},
the function $\sin(\lambda)$ is strongly realizable,
since $\sin(\lambda)=\sqrt{1-\cos(\lambda)^2}$.

\subsection{Algebraic functions and stochastic grammars.}
It is well known that unambiguous context-free grammars are associated with generating functions
that are algebraic: see~\cite{FlSe09} for background (the Chomsky--Sch\"utzenberger Theorem).

\begin{definition} A \emph{binary stochastic grammar} (a ``bistoch'') is a context-free grammar 
whose terminal alphabet is binary, conventionally $\{H,T\}$,
where each production is of the form
\begin{equation}\label{cf0}
\cal X \longrightarrow H {\frak m} + T {\frak n},
\end{equation}
with ${\frak m},{\frak n}$ that are monomials in the non-terminal symbols. It is assumed that each non-terminal is
the left hand side of at most one production.
\end{definition}
Let $G$, with axiom~$\cal S$, be a bistoch. We let ${\bf L}[G;\cal S]$ be the language 
associated with~$\cal S$. By the classical theory of formal languages and automata,
this language can be recognized by a pushdown (stack) automata.
The constraint that there is a single production for each non-terminal on the left
means that the automaton corresponding to a bistoch is \emph{deterministic}.
(It is then a simple matter to come up with a recursive procedure that
parses a word in $\{H,T\}^\star$ according to a non-terminal symbols~$S$.)
In order to avoid trivialities, we assume that all non-terminals are ``useful'',
meaning that each produces at least one word of $\{H,T\}^\star$.
For instance, the one-terminal grammar
$ 
\cal Y = H \cal Y \cal Y \cal Y  + T
$ 
generates all \L ukasiewicz codes of ternary trees~\cite[\S I.5.3]{FlSe09} 
and is closely related to the construction of Note~\ref{pda-note}. 

Next, we introduce the
\emph{ordinary generating function} (or \ogf) of $G$ and~$S$,
\[
S(z):=\sum_{w\in {\bf L}[G;\cal S]} z^{|w|} =\sum_{n\ge0} S_n z^n,
\]
with $S_n$ the number of words of length~$n$ in ${\bf L}[G;\cal S]$.
The deterministic character of a bistoch grammar implies that 
the  {\ogf}s are bound by a system of equations
(one for each nonterminal): from~\eqref{cf0}, we have
$ 
X(z)=z \wh {\frak m} + z \wh {\frak n},
$ 
where $\wh{\frak m}, \wh{\frak n}$ are monomials in the {\ogf}s 
corresponding to the non-terminals of~${\frak m},{\frak n}$; see~\cite[\S I.5.4]{FlSe09}. 
For instance, 
in the ternary tree case: 
$Y=z+zY^3$.

Thus, any {\ogf}~$y$ arising from a bistoch is a component 
of a system of polynomial equations, hence,
an \emph{algebraic function}. By elimination, the system
reduces to a single 
equation $P(z,y)=0$.
We obtain, with a simple proof, a result of Mossel-Peres~\cite[Th.~1.2]{MoPe05}:

\begin{figure}\small

{\rm \begin{tabbing}\rm
\= XX \= XX\= XX \= XXX\=\kill
\> $\{$\,{\bf let} $N:=\ggeo(\lambda)$; \\
\> {\bf draw} $w:=X_1X_2\cdots X_{N}$ {with} $\Pr(X_j=H)=\Pr(X_j=T)=\frac12$;\\
\> \iff~$w\in{\bf L}(G;\cal S)$ \thenn~return(1) \elsee~return(0)\,$\}$.
\end{tabbing}}%
\vspace*{-2.5truemm}
\caption{\label{bistoch-fig}\small
The algebraic construction associated to the pushdown 
automaton arising from a bistoch grammar.}
\end{figure}

\begin{theorem}[\cite{MoPe05}]
To each bistoch grammar~$G$ and non-terminal~$\cal S$, there corresponds 
a construction (Fig.~\ref{bistoch-fig}),
which can be implemented by a deterministic pushdown automaton and calls to a~$\gber(\lambda)$
and is of type
$ 
\gber(\lambda)~\longrightarrow ~\gber\left(S\left(\frac{\lambda}{2}\right)\right),
$ 
where $S(z)$ is the algebraic function canonically associated with 
the grammar~$G$ and non-terminal~$S$.
\end{theorem}

\begin{note} \emph{Stochastic grammars and positive algebraic functions.}
First, we observe that another way to describe the process is by means of
a stochastic grammar with production rules
$
\cal X ~ \longrightarrow~ \hbox{$\frac12$} {\frak m} + \hbox{$\frac 12$} {\frak n},
$
where each possibility is weighted by its probability ($1/2$).
Then fixing~$N=n$ amounts to
conditioning on the size~$n$ of the resulting object.
This bears a superficial resemblance to branching processes,
upon conditioning on the size of the total progeny, itself assumed to be finite.
(The branching process may well be supercritical, as in the ternary tree case.) 

The algebraic generating functions that may arise from such grammars and 
positive systems of equations have been widely studied. Regarding coefficients
and singularities, we refer to the discussion of the 
Drmota--Lalley--Woods Theorem in~\cite[pp.~482--493]{FlSe09}.
Regarding the values of the generating functions, we mention the
studies by Kiefer, Luttenberger, and Esparza~\cite{KiLuEs07} and
by Pivoteau, Salvy, and Soria~\cite{PiSaSo08}. The former is motivated by 
the probabilistic verification of recursive Markov processes,
the latter by the efficient implementation of Boltzmann samplers.

It is not (yet) known whether a function as simple as
$(1-\lambda)^{-1/3}$ is realizable by a stochastic context-free grammar
or, equivalently, a deterministic pushdown automaton. (We conjecture 
a negative answer, as it seems that only square-root and iterated square-root singularities are 
possible.)\par 
\end{note}

\subsection{Ramanujan, hypergeometrics, and a Buffon computer for $1/\pi$.}
A subtitle 
might be: \emph{What to 
do if you want to perform Buffon's experiment but don't have needles, just coins?}
The identity 
\[
\frac{1}{\pi}=\sum_{n=0}^\infty \binom{2n}{n}^3 \frac{6n+1}{2^{8n+4}},
\]
due to Ramanujan (see~\cite{Guillera06} for related material), 
lends itself somewhat miraculously to evaluation by a 
simple Buffon computer.
The following simple experiment (the probabilistic procedure \emph{Rama})
succeeds (returns the value~1) with probability exactly $1/\pi$.
It thus constitutes a discrete analogue of Buffon's original,
one with only three counters ($T$, a copy of~$T$, and~$\Delta$).

\begin{tabbing}\small
n\= XXl \= l \= l \= XXX \= \kill
\> {\sf procedure} Rama();~~~\{\emph{returns~1 with probability $1/\pi$}\} \\
\> {\bf S1.} \> {\sf let} $T$ := $X_1+X_2$, {\sf where} $X_1,X_2\in\operatorname{Geom}(\frac14)$; 
\\ 
\> {\bf S2.} \> {\sf with} probability $\frac59$ {\sf do} $T$ := $T+1$; \\
\> {\bf S3.} \> {\sf for} $j=1,2,3$ {\sf do} \\
\> {\bf S4.} \> \ {\sf draw}  a sequence of $2T$ coin flips; \\
\> \> \  {\sf if} ($\Delta\equiv{}$ \# Heads${}-{}$\# Tails)\,${}\not=0$ {\sf then} return(0); \\
\> {\bf S5.} \> return(1).
\end{tabbing}


\section{\bf A Buffon integrator}\label{int-sec}


Our purpose here is to develop, given a construction of type
$\gber(\lambda) ~\longrightarrow~ \gber(\phi(\lambda))$, a generator for the function 
\begin{equation}\label{int1}
\Phi(\lambda)=\frac{1}{\lambda} \int_0^\lambda \phi(w)\, dw.
\end{equation}
An immediate consequence will be a generator for $\lambda\Phi(\lambda)$; that is,
an \emph{``integrator''}. 

To start with, we discuss a purely discrete implementation of
$
\gber(\lambda)~\longrightarrow~ \gber(U\lambda)$, with $U\in[0,1]$, uniformly,
where multiple invocations of $\gber(\lambda)$ must involve the same value of~$U$.
Conceptually, it suffices to draw $U$ as an infinite sequence of flips,
then make use of this~$U$ to get a $\gber(U)$ and then appeal to the conjunction (product)
construction to get a $\gber(\lambda U)$ as $ \gber(U)\cdot \gber(\lambda)$. 
%
\begin{figure}\small

\begin{tabular}{cc}
\qquad\begin{tabular}{c}\setlength{\unitlength}{0.75truecm}
\begin{picture}(7.7,5)
\put(-1,2){$U={}$}
\put(0,0){\Img{3}{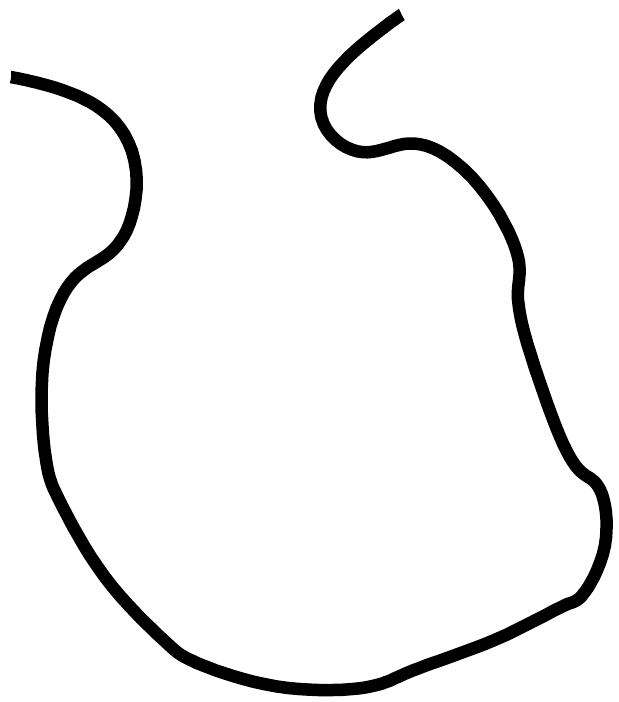}}
\put(0.8,1.3){$\bf 1_1$}
\put(2.7,0.8){$\bf 2_0$}
\put(1.4,2.9){$\bf 3_0$}
\put(2.0,1.7){$\bf 5_0$}
\put(2.5,2.3){$\bf 6_1$}
\thicklines
\put(0.9,1.4){\circle{0.7}}
\put(2.8,0.9){\circle{0.7}}
\put(1.5,3.0){\circle{0.7}}
\put(2.1,1.8){\circle{0.7}}
\put(2.6,2.4){\circle{0.7}}
\put(4,2){${}={}$}
\put(6,2.2){
$\begin{array}{|c|}
\hline
?\\?\\?\\ \bf 1\\ \bf 0 \\ ? \\ \bf 0\\ \bf 0 \\ \bf 1 \\ \hline\hline
\end{array}$}
\put(5.1,2.25){
$\begin{array}{c}
9~:\\ 8~:\\ 7~:\\ 6~:\\ 5~:\\ 4~:\\ 3~:\\ 2~:\\ 1~:
\end{array}$}
\put(6.9,3.3){$\Longleftarrow~~J$}
\end{picture}\end{tabular}
&
%
%
\end{tabular}
\vspace*{-0.45truecm}

\caption{\label{X-fig}\small
The ``geometric-bag''  procedure~$bag(U)$:
 two graphic representations of a state (the pairs index--values and
 a partly filled register).
}
\end{figure}
To implement this idea, it suffices to resort to \emph{lazy evaluation}.
One may think of~$U$ as a potentially infinite vector $(\up_1,\up_2,\ldots)$,
where $\up_j$ represents the $j$th bit of~$U$. Only a finite section
of the vector is used at any given time and the $\up_j$ are initially undefined
(or unspecified). Remember that a $\gber(U)$ is simply obtained by fetching the bit of~$U$
that is of order $J$, where $J\in1+\geo(\frac12)$; cf Eq.~\eqref{bern0}. 
In our relaxed lazy context, whenever such a bit~$\up_j$ is fetched, we first
examine whether it has already been assigned a $\{0,1\}$--value; if so,
we return this value; if not, we first perform a flip, assign it to~$\up_j$,
and return the corresponding value: see Fig.~\ref{X-fig}. (The implementation is 
obvious: one can maintain an association list of the already ``known'' indices and values,
and update it, as the need arises; or keep a boolean vector
of not yet assigned values; or encode a yet unassigned position
by a special symbol, such as `?' or `$-1$'. See the Appendix for 
a simple implementation.)

Assume that $\phi(\lambda)$ is realized by a Buffon machine 
that calls a Bernoulli generator 
$\gber(\lambda)$. If we replace $\gber(\lambda)$ by 
$\gber(\lambda U)$, as described in the previous paragraph, 
we obtain a Bernoulli generator whose parameter is $\phi(\lambda U)$, where $U$ is uniform over $[0,1]$. This is equivalent to a Bernoulli generator whose parameter is
$
\int_0^1 \phi(\lambda u)\, du = \Phi(\lambda)$,
with~$\Phi(\lambda)$ as in~\eqref{int1}. 


\begin{theorem}\label{int-thm}
Let $\phi(\lambda)$ be realizable by a Buffon machine~$\cal M$.
Then the function
$\Phi(\lambda)=
\frac{1}{\lambda} \int_0^\lambda \phi(w)\, dw$ is realizable
by addition of a geometric bag to $\cal M$.
In particular, if $\phi(\lambda)$ is realizable,
then its integral taken starting from~$0$ is also realizable.
\end{theorem}
This result paves the way to a large number of derived constructions.
For instance, starting from the even-parity construction
of~\S\ref{frame-sec}, we obtain 
$ 
\Phi_0(\lambda):=\frac{1}{\lambda}\int_0^\lambda \frac{1}{1+w}\, dw =\frac{1}{\lambda} \log(1+\lambda),
$ 
hence, by product, a construction for $\log(1+\lambda)$.
 When we now combine the parity construction with 
``squaring'', where a $\gber(p)$ is replaced by the product
$\gber(p)\cdot\gber(p)$, we obtain
$ 
\Phi_1(\lambda):=\frac{1}{\lambda}\int_0^\lambda \frac{dw}{1+w^2}=\frac{1}{\lambda}\arctan(\lambda),
$ 
hence also $\arctan(\lambda)$.
When use is made of the exponential (Poisson) construction $\lambda\mapsto e^{-\lambda}$,
one obtains (by squaring and after multiplications) a construction for
$ 
\Phi_2(\lambda):=\int_0^\lambda e^{-w^2/2}\, dw,
$ 
so that the error function (``erf'') is also realizable. Finally,
the square-root construction 
combined with parity and integration provides
$ 
\Phi_3(\lambda):=\int_0^\lambda \frac{\sqrt{1-w^2}}{1+w}\, dw = -1+\sqrt{1-\lambda^2}+\arcsin(\lambda),
$ 
out of which we can construct $\frac12\arcsin(\lambda)$.
In summary:
\begin{theorem}\label{invtrig-thm}
The following functions are strongly realizable $(0\le x<1)$:
\[
\log(1+x),~~\arctan(x),~~\frac12\arcsin(x),~~\int_0^x e^{-w^2/2}\, dw.
\]
The first two only require one bag; the third requires a bag and a stack;
the fourth can be implemented with a string register and bag.
\end{theorem}
\parag{Buffon machines for $\pi$.}
The fact that
$\ds 
\Phi_1(1)=\arctan(1)=\frac{\pi}{4}.
$ 
yields a Buffon computer for $\pi/4$. 
There are further simplifications due to the fact that $\gber(1)$ is trivial:
this computer then only makes use of the~$U$ vector.
Given its extreme simplicity, we  can even list the complete code of this
 Madhava--Gregory--Leibniz (MGL) generator for $\pi/4$:
\begin{small}
\begin{verbatim}
    MGL:=proc()    do 
          if bag(U)=0 then return(1) fi;  
          if bag(U)=0 then return(1) fi;
          if bag(U)=0 then return(0) fi;  
          if bag(U)=0 then return(0) fi; od; end.
\end{verbatim}
\end{small}%
%
\noindent
The Buffon computer based on $\arctan(1)$ works fine for small simulations.
For instance, based on 10,000 experiments, we infer the 
approximate value $\pi/4\approx {\bf 0.78}76$, whereas $\pi/4\doteq {\bf 0.78}539$,
with a mean number of flips per experiment  about~$27$.
However, 
values of $U$ very close to~$1$ are occasionally
generated (the more so, as the number of simulation increases).
Accordingly, 
 \emph{the expected number of flips is infinite},
a fact to be attributed to slow convergence in
the Madhava--Gregory--Leibniz series,
$
\frac{\pi}{4} = \frac{1}{1}-\frac{1}{3}+\frac{1}{5}-\frac{1}{7}+\cdots\,.
$ 

%

\begin{figure}

\hbox{\ \Img{3.85}{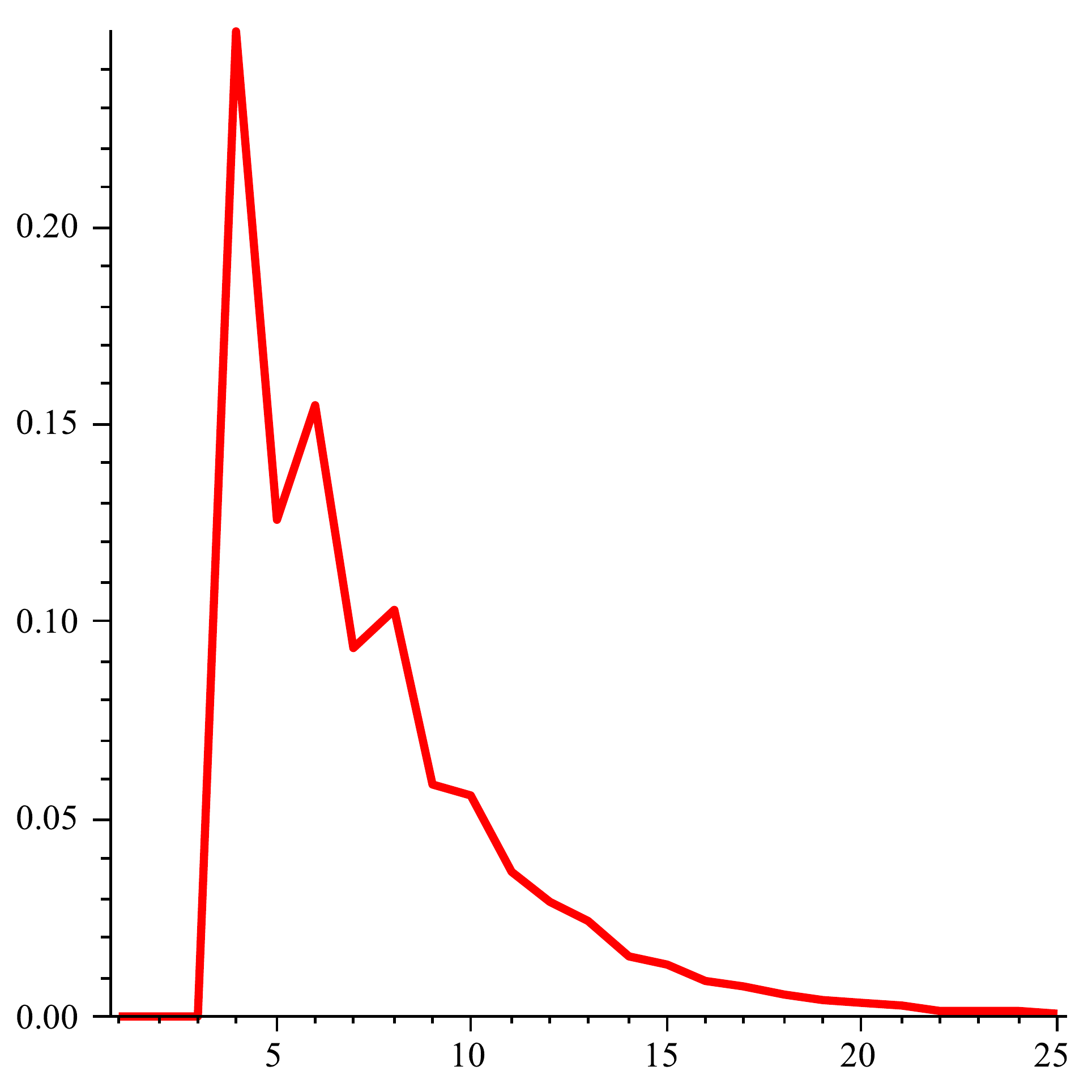}~
\Img{3.85}{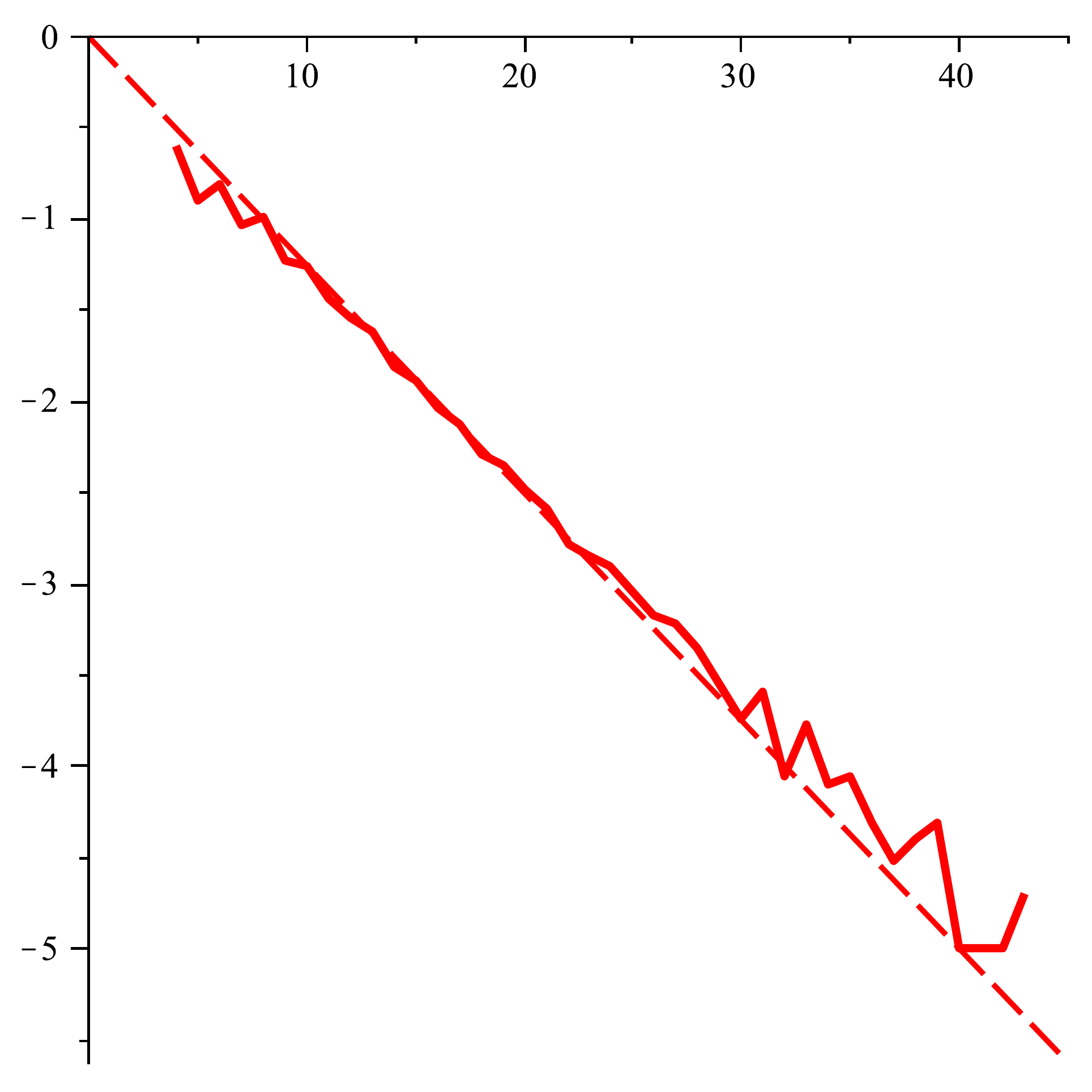}}

\vspace*{-0.5truemm}
\caption{\label{machinc-fig}\small
The distribution of costs of the Machin machine~\eqref{machin23b}.
\emph{Left}:  histogram.
\emph{Right}: decimal logarithms of the probabilities,
compared to $\log_{10}(10^{-k/8})$ (dashed line).
}
\end{figure}

%

The next idea is to consider formulae of a kind made well-known by Machin,
who appealed to arc-tangent addition formulae
in order to reach the record computation of 100 digits of~$\pi$ in 1706.
For our purposes, a formula without negative signs is needed,
the simplest  of which,
\begin{equation}\label{machin23}
\frac{\pi}{4}=\arctan\left(\frac{1}{2}\right)+\arctan\left(\frac{1}{3}\right),
\end{equation}
being especially suitable for a short program
 is easily compiled \emph{in silico} under the form
\begin{equation}\label{machin23b}
\qquad \frac{\pi}{4}=\frac{1}{2}\left[2\arctan\left(\frac12\right)
+\frac23\cdot 3\arctan\left(\frac13\right)\right].
\end{equation}
(This last form only uses the realizable functions $\lambda^{-1}\arctan(\lambda)$,
$2\lambda/3$ and the binary operation $\frac12[p+q]$.)

With $10^6$ simulations, we obtained an estimate $\pi/4\approx {\bf 0.785}98$,
to be compared to the exact value~$\pi/4={\bf0.785}39\cdots$; that 
is, an error of about~$6\cdot 10^{-4}$, well within
normal statistical fluctuations. The empirically measured average number of flips 
per generation of this Machin-like $\gber(\pi/4)$ 
turned out to be about $6.45$ coin flips. 
Fig.~\ref{machinc-fig} furthermore displays the empirical distribution of the number of
coin flips, based on another campaign of~$10^{5}$ simulations.
The distribution of costs appears to have exponential tails matching 
fairly well the approximate formula~$\Pr(C=k)\approx 10^{-k/8}$.
The complete code for a version 
of this generator, which produces $\frac\pi8$, is given in the Appendix.

%

Yet an alternative construction is based on the arcsine and $\Phi_3$. 
%
Many variations are possible, related to multiple or iterated integrals (use several bags).

\section{\bf Experiments}\label{concl-sec}
\begin{figure}[t]\small

\begin{center}
\Img{8.5}{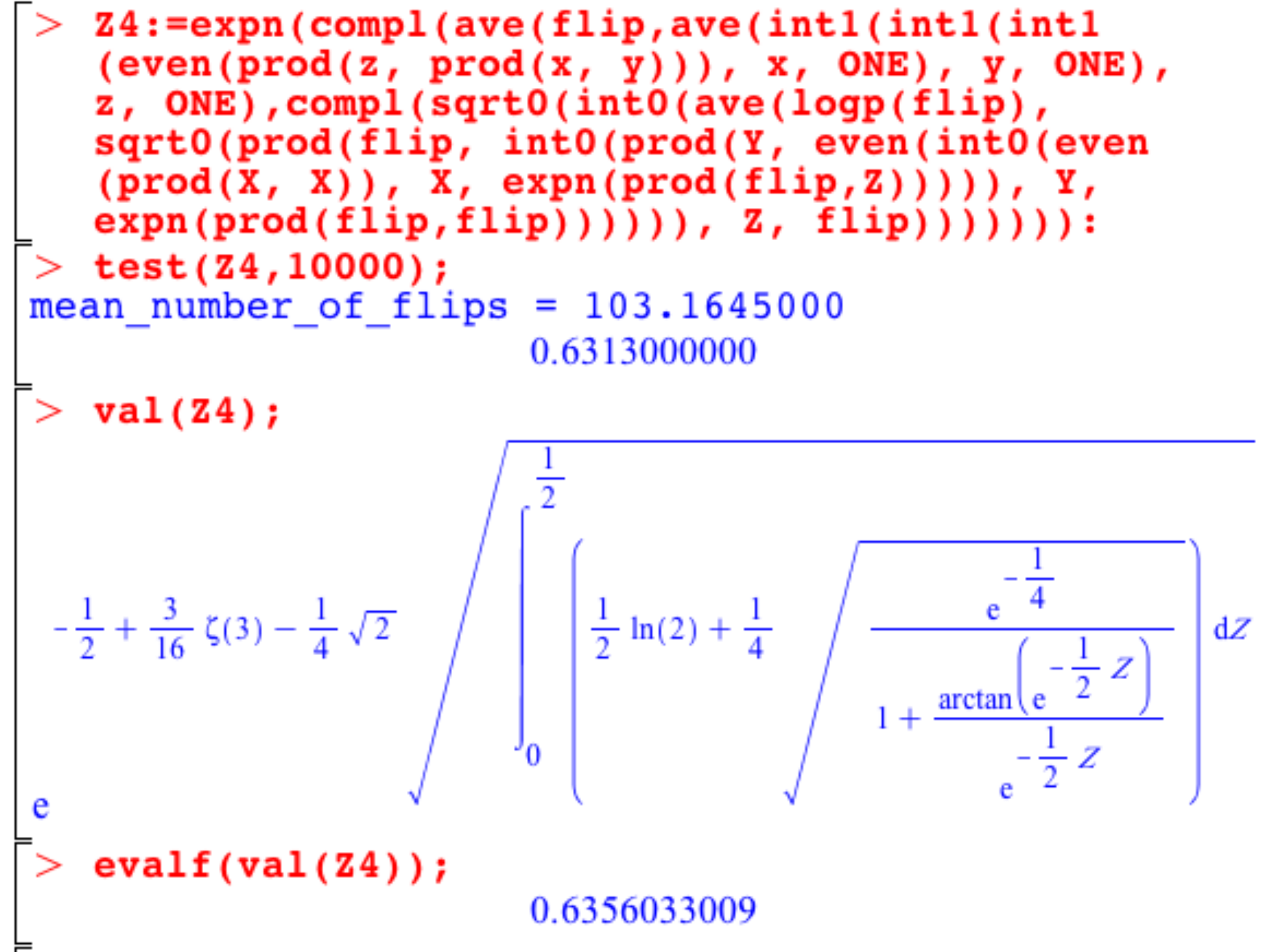}
\end{center}
\vspace{-0.25truecm}
\caption{\label{screen-fig}\small Screen copy of a {\sc Maple} session fragment showing:
$(i)$ the symbolic description of a generator;
$(ii)$~a simulation of $10^4$ executions having a proportion of successes equal to ${\bf 0.63}130$,
with a mean number of flips close to 103;
$(iii)$~the actual symbolic value of the probability of success
and its numerical evaluation~${\bf 0.63}560\cdots$\,.
}
\end{figure}

We have built a complete prototype implementation under the {\sc Maple} symbolic manipulation system,
in order to test and validate the ideas expounded above; see Fig.~\ref{screen-fig}. 
A generator, such as~$Z4$ 
of Fig.~\ref{screen-fig}, is specified as a composition of basic constructions, such
as $f\mapsto \exp(-f)$ [{\tt expn}], $f\mapsto\sqrt{f}$ [{\tt sqrt0}], 
$f\mapsto \int f$ [{\tt int1}], and so on. 
An interpreter using as source of randomness the built-in function {\tt random} 
then takes this description and produces a $\{0,1\}$ result;
this interpreter, which is comprised of barely 60 lines, contains from one to about a dozen
instructions for each construction. In accordance with the conditions of our contract,
only simple register manipulations are ever used, so that a transcription in~C or Java
would be of a roughly comparable size.

We see here that even a complicated constant such as the ``value'' of the
probability associated with~$Z4$, 
\[\hbox{$
e$\footnotesize${}^{-\frac12+\frac3{16}\,\zeta  \left( 3 \right) -\frac14\,\sqrt {2}\sqrt {{\ds\int _{0}
^{\frac12}}\!\frac12\ln   2  +\frac14\,\sqrt {{e^{-\frac14}} \left( 1+{
\frac {\operatorname{atan} \left( {e^{-Z/2}} \right) }{{e^{-Z/2}}}} \right)
^{-1}}\!{dZ}}}$},
\]
 is effectively simulated with an error of~$4\, 10^{-3}$, which is once more
consistent with normal statistical fluctuations. For such a complicated constant, the
observed mean number of flips is a little above 100. 
Note that the quantity $\zeta(3)$ is produced (as well as retrieved automatically by
{\sc Maple}'s symbolic engine!) from Beuker's triple integral:
$
\frac{7}{8}\zeta(3)=\int_0^1\int_0^1\int_0^1
\frac{1}{1+xyz}\, dx\, dy\, dz.
$
(On batches of $10^5$ experiments, that quantity alone only consumed an average of $6.5$ coin flips,
whereas the analogous $\frac{31}{32}\zeta(5)$  required barely~6 coin flips on average.)

Note that the implementation is, by nature, freely extendible.
Thus, given integration and our basic primitives (e.g., even(f)${}\equiv\frac{1}{1+f}$),
we readily program an arc-tangent as a one-liner, 
\begin{equation}\label{atan}
\arctan(f) = f\cdot \left[\frac{1}{f} \int_0^f \frac{dx}{1+x^2}\right],
\end{equation}
and similarly for sine, arc-sine, ~$\operatorname{erf}$, etc, with the symbolic engine automatically providing
the symbolic and numerical values, as a validation.

Here is finally a table recapitulating nine ways of building Buffon machines for 
$\pi$-related
constants, with, for each of the methods,  the value, and 
empirical average of the number of coin flips, as observed
over $10^4$ simulations:
\begin{displaymath}\label{pitab}\renewcommand{\tabcolsep}{2pt}
\renewcommand{\arraystretch}{1.4}
\hbox{\footnotesize\begin{tabular}{c|c|ccc|cc|cc}
\hline\hline
$\operatorname{Li}_2(\frac12)$ & Rama & 
\multicolumn{3}{c|}{arcsin $[1;\frac{1}{\sqrt{2}};\frac12]$} & \multicolumn{2}{c|}{arctan $[\frac12+\frac13; 1$]} &
$\zeta(4)$ & $\zeta(2)$\\
\hline
$\ds\frac{\pi^2}{24}$ & $\ds\frac{1}{\pi}$ & $\ds\frac{\pi}{4}$ & $\ds\frac{\pi}{4}$ & $\ds\frac{\pi}{12}$ 
 & $\ds\frac{\pi}{4}$ & $\ds\frac{\pi}{4}$ &  $\ds\frac{7\pi^4}{720}$ & $\ds\frac{\pi^2}{12}$ \\[2truemm]
7.9 & 10.8 & 76.5 ($\infty$) & 16.2 & 4.9 & 6.5 & 26.7 $(\infty)$ & 6.2 & 7.2.\\
\hline\hline
\end{tabular}
}
\end{displaymath}
(The tag ``$\infty$'' means that the expected cost of the simulation is infinite---a weak realization.)

\section{Conclusion}


As we pointed out in the introduction, every computable number
can be simulated by a machine, but one that, in general, will violate
our charter of simplicity (as measured, typically, by program size).
Numbers accessible to our framework
seem not to include Euler's constant $\gamma\doteq 0.57721$, and we must
leave it as an open problem to come up with a ``natural'' experiment,
whose probability of success is~$\gamma$. Perhaps the difficulty
of the problem lies in the absence of a simple ``positive'' expression
that could be compiled into a correspondingly simple Buffon generator.
By contrast, exotic numbers, such as
$\pi^{-1/\pi}$ or $e^{-\sin(1/\sqrt{7})}$ are easily simulated\ldots

On another  note, we have  not considered the  generation of
\emph{continuous}  random variables~$X$,  specified by a  distribution
function $F(x)=\Pr(X\le x)$.   Von Neumann's original algorithm for an
exponential variate belongs to this  paradigm.  In this case, the bits
of~$X$ are  obtained by a  short computation of~$O(1)$ initial
bits,  continued by  the    production   of  an  infinite   flow    of
\emph{uniform} random  bits.   This  theme is  thoroughly  explored by
Knuth and Yao in~\cite{KnYa76}.  It would be of obvious interest to be
able to hybridize the  von Neumann-Knuth-Yao generators of  continuous
distributions with  our Buffon  computers for  discrete distributions.
Interestingly, the fact that the Gaussian error function, albeit
restricted to the interval~$(0,1)$,  is realizable by Buffon  machines
suggests the possibility of a totally  discrete generator for a
(standard) \emph{normal} variate.

The present work  is
excellently   summarized   by  Keane   and   O'Brien's  vivid expression of   
\emph{``Bernoulli      factory''}~\cite{KeOB94}. 
It was initially approached with  a purely theoretical goal. 
It then came  as   a surprise that  \emph{a  priori}  stringent theoretical
requirements---those of perfect generation, discreteness of the random
source,     and  simplicity   of   the    mechanisms---could lead   to
computationally efficient algorithms. We have  indeed seen many cases,
from Bernoulli to    logarithmic and Poisson generators,   where short
programs and the execution of just a few dozen instructions suffice!


\smallskip

\begin{small}
\noindent{\bf Acknowledgements.} This work was supported by 
the French ANR Projects {\sc Gamma}, {\sc Boole}, and {\sc Magnum}.
The authors are grateful to the referees of SODA-2011 for their 
perceptive and encouraging comments.\par 
\end{small}


\bibliographystyle{acm}
\bibliography{algo}

\begin{appendix}\small

\section*{Appendix: A complete  Buffon machine for $\frac{\pi}{8}$}

Here is the complete pseudo-code 
(in fact an executable {\sc Maple} code), cf procedure {\tt Pi8} below,
for a ${\pi}/{8}$ experiment, 
based on Eq.~\eqref{machin23}  
and not using any high-level primitive.
It exemplifies many constructions seen in the text. The translation to
various low level languages, such as~{\sf C}, should be immediate, 
possibly up to inlining of code or macro expansions, 
in case procedures cannot be passed as arguments
of other procedures. The expected number of coin flips per experiment
is about 4.92.


\medskip

\hrule

\medskip
\noindent
The \emph{flip} procedure (returns a pseudo-random bit):
\begin{verbatim}
flip:=proc() if rand()/10.^12<1/2 
   then return(1) else return(0) fi; end:
\end{verbatim}
A $\ggeo(\frac12)$ returns a geometric of parameter $\frac12$; cf Eq.~\eqref{geom}:
\begin{verbatim}
1.  geo_half:=proc () local j;  j := 0;
2.    while flip() = 0 do j:=j+1 od; return j end:
\end{verbatim}
A $\gber(\frac13)$ returns a Bernoulli of parameter $\frac13$; cf Eq.~\eqref{bern0}:
\begin{verbatim}
3.  bern_third:=proc() local a,b; 
4.    do a:=flip(); b:=flip(); 
5.      if not ((a=1) and (b=1)) then break fi; od; 
6.    if (a=0) and (b=0) then return(1) 
7.        else return(0) fi; end:
\end{verbatim}
\emph{Bags}. Initialization and result of a comparison with a random $U\in[0,1]$;
cf Fig.~\ref{X-fig} and~\S\ref{int-sec}:
\begin{verbatim}
    INFINITY:=50: 
8.  init_bagU:=proc() local j; global U; 
9.    for j from 1 to INFINITY do U[j]:=-1 od; end:
10. bagU:=proc() local k;global U; k:=1+geo_half(); 
11.  if U[k]=-1 then U[k]:=flip() fi; return(U[k]);
                                               end:
\end{verbatim}
\vspace*{-3truemm}
(To obtain a \emph{perfect} generator, dynamically increase {\tt INFINITY},
if ever needed---the probability is $<10^{-15}$.)
\smallskip

The {\tt EVENP} construction takes $f\in\gber(\lambda)$ and produces a
 $\gber(\frac{1}{1+\lambda})$; cf Eq.~\eqref{evenp}:
\begin{verbatim}
12. EVENP:=proc(f) do if f()=0 then return(1); fi;
13.        if f()=0 then return(0) fi; od; end:
\end{verbatim}
The main {\tt atan} construction, based on bags implementing integration,
takes an $f\in \gber(\lambda)$
and produces a $\gber(\arctan(\lambda))$;
the auxiliary procedure $g()$ builds a $\gber(\lambda^2 U^2)$,
with $U\in[0,1]$ random; cf \S\ref{int-sec} and Eq.~\eqref{atan}: 
\begin{verbatim}
14. ATAN:=proc(f) local g; 
15.  if f()=0 then return(0) fi;
16.  init_bagU();
17.  g:=proc() if bagU()=1 then if f()=1 then 
18.     if bagU()=1 then return(f()) fi; fi; fi; 
19.     return(0); end;
20.  EVENP(g); end:
\end{verbatim}
The {\tt Pi8} procedure is a $\gber(\frac{\pi}{8})$
based on Machin's arc tangent formula
(the arithmetic mean of
$\arctan(1/2)$ and $\arctan(1/3)$ is taken); cf Eq.~\eqref{machin23}:
\begin{verbatim}
21. Pi8:=proc() if flip()=0 then ATAN(flip) 
22.             else ATAN(bern_third) fi end:.
\end{verbatim}

\hrule

\end{appendix}

\end{twocolumn}
\end{document}